\documentclass[12pt,reqno]{article}

\usepackage[usenames]{color}
\usepackage{amssymb}
\usepackage{amsmath}
\usepackage{amsthm}
\usepackage{amsfonts}
\usepackage{amscd}
\usepackage{graphicx}
\usepackage{chngcntr}

\usepackage[colorlinks=true,
linkcolor=webgreen,
filecolor=webbrown,
citecolor=webgreen]{hyperref}

\definecolor{webgreen}{rgb}{0,.5,0}
\definecolor{webbrown}{rgb}{.6,0,0}

\usepackage{color}
\usepackage{fullpage}
\usepackage{float}

\usepackage{graphics}
\usepackage{latexsym}
\usepackage{epsf}
\usepackage{breakurl}
\usepackage[english]{babel}

\numberwithin{equation}{section}
\numberwithin{figure}{section}

\usepackage{tikz}

\usepackage{hyperref}

\makeatother

\usepackage{babel}

\theoremstyle{plain}
\newtheorem{theorem}{Theorem}
\newtheorem{corollary}[theorem]{Corollary}

\newtheorem{proposition}[theorem]{Proposition}
\theoremstyle{definition}
\newtheorem{definition}[theorem]{Definition}
\newtheorem{example}[theorem]{Example}
\theoremstyle{remark}
\newtheorem{remark}[theorem]{Remark}

\begin{document}

\newcommand{\seqnum}[1]{\href{https://oeis.org/#1}{\rm \underline{#1}}}

\title{Number of Partitions of Modular Integers\\
(with an Appendix by P.\ Deligne)} 

\date{}

\author{David Broadhurst and Xavier Roulleau} 

\maketitle

\begin{abstract}
For integers $n,k,s$, we give  a formula for the number $T(n,k,s)$ 
of order $k$ subsets of the ring $\mathbb{Z}/n\mathbb{Z}$ whose
sum of elements is $s$ modulo $n$. To do so, we describe   
explicitly a sequence of matrices $M(k)$, for positive integers $k$,  
such that the size of $M(k)$ is the number of divisors of $k$, and
for two coprime integers $k_{1},k_{2}$, the matrix $M(k_{1}k_{2})$
is the Kronecker product of $M(k_{1})$ and $M(k_{2})$. For $s=0, 1, 2$,  
and for $s=k/2$ when $k$ is even, the sequences $T(n,k,s)$ are related  
to the number of necklaces with $k$ black beads and $n-k$ white  
beads, and to Lyndon words. This work begins with empirical
determinations of $M(k)$ up to $k=10000$, from which we infer a closed formula
that encompasses many entries in the Encyclopedia of Integer Sequences.
Its proof comes from work on Ramanujan sums, by Ramanathan,  with
a generalization to wider problems linked to representation theory and recently described by Deligne.
\end{abstract}

\begin{center}
\noindent David Broadhurst \\ 
The Open University \\
School of Physical Sciences\\
Milton Keynes MK7 6AA\\ 
UK\\ 
\href{mailto:david.broadhurst@open.ac.uk}{david.broadhurst@open.ac.uk}

\vspace{0.3cm} 
\noindent Xavier Roulleau\\
Universit\'e d'Angers \\
CNRS, LAREMA, SFR MATHSTIC \\
F-49000 Angers\\ 
France \\ 
\href{mailto:xavier.roulleau@univ-angers.fr}{xavier.roulleau@univ-angers.fr}
\end{center}

\tableofcontents{}

\section{Introduction}

For integers $n$, $k$, and $s$, with $n\ge k>0$, let 
\[
S(n,k,s)
\]
denote the set of order $k$ subsets of the set of integers $\{1,\dots n\}$  
(which can also be seen as the ring $\mathbb{Z}/n\mathbb{Z}$) whose 
sum of elements is congruent to $s$ modulo $n$. We denote by 
\[
T(n,k,s)=|S(n,k,s)|
\]
the cardinality of $S(n,k,s)$. The aim of this paper is to describe
a formula for $T(n,k,s)$ and to relate it to entries
in the On-line Encyclopedia of Integer Sequences (OEIS) founded by Sloane \cite{Sloane}. 

We first introduce some definitions and notation.
For $q$ a power of a prime integer and for positive integers $u,v$ dividing  
$q$, let us define 
\begin{equation}
M(q)_{u,v}=\begin{cases}
\mu(v/u)u, & \text{if \ensuremath{v>u}};\\ 
\varphi(v), & \text{otherwise}  
\end{cases}\label{Muv}
\end{equation}
where $\mu$ is the M\"obius function and $\varphi$ is Euler's totient
function. 
\begin{theorem}
\label{thm:MainIntro}For integers $n$, $k$ and $s$, with $n\ge k>0$,
one has 
\[
T(n,k,s)=\frac{1}{n}\sum_{d\mid\gcd(n,k)}(-1)^{k-k/d}M(k)_{\gcd(k,s),d}\binom{n/d}{k/d}
\]
where $M(k)_{t,d}$ with $t\mid k$ and $d\mid k$ is the integer computed
by factorizing $k$ into prime powers $q=p^{e}$, with $e$ the valuation  
of $k$ at prime $p|k$, to give 
\begin{equation}
M(k)_{t,d}=\prod_{q\mid k}M(q)_{\gcd(q,t),\gcd(q,d)}.\label{Mtd}
\end{equation}
\end{theorem}

It is easy to show that $T(n,k,s)$ is equal to $T(n,k,\gcd(k,s))$.  
The integers $M(k)_{t,d}$ for $t\mid k,\,d\mid k$ form a square matrix $M(k)$
indexed by the divisors of $k$, thus of size the number of divisors
of $k$. These matrices have remarkable properties, which we describe
in Theorem \ref{thm:Merveilles}, and are perhaps the main interest
of the paper. In particular, we prove the following theorem. 
\begin{theorem}\label{Theorem2} 
Let $k,k'$ be two coprime positive integers. The
matrix $M(kk')$ is the Kronecker product of $M(k)$ and $M(k')$.
\end{theorem}

As explained in the text, since the matrices are indexed by divisors
of coprime integers, the Kronecker product is commutative: 
$M(k)\otimes M(k')=M(k')\otimes M(k)=M(kk')$.

Now, let us provide some further motivations
regarding the numbers $T(n,k,s)$
and discuss examples of related sequences in the OEIS.
The case $k=3$ and $s=0$ has a natural geometric interpretation,
as follows. Let $E\hookrightarrow\mathbb{P}^{2}$ be an elliptic curve: a
smooth plane cubic curve with a distinguished flex point $O$, so
that there is an additive group law on the points of $E$, with neutral
element $O$. This group law is constructed in such a way that  
three points $p_{1},p_{2},p_{3}$ on $E$ lie on a line if and only
if $p_{1}+p_{2}+p_{3}=O$. Let $p$ be a point on $E$ of order exactly
$n$, so that $G_{n}=(jp)_{j\in\mathbb{Z}/n\mathbb{Z}}$ is isomorphic
to $\mathbb{Z}/n\mathbb{Z}$. The set of lines in the plane that contain
exactly $3$ points of $G_{n}$ is in bijection with $S(n,3,0)$,
and the number of such lines, $T(n,3,0)$, as $n$ varies, is sequence
\seqnum{A007997} in the OEIS

Then, there is a natural generalization for conics and order six  
elements of $\mathbb{Z}/n\mathbb{Z}$: a conic (possibly degenerate,
i.e., the union of two lines) cuts the cubic curve $E$ at six points
(with multiplicities) if and only if their sum is $O$. Thus $T(n,6,0)$
is also the number of conics containing six points in $G_{n}$. Moreover, sequence  
\seqnum{A379920} is the number of \emph{irreducible} conics containing six 
points in $G_{n}$ for $n\geq9$. 

We did not find the sequence $(T(n,6,0))_{n\geq 6}$ in the OEIS,   
but have now added it as \seqnum{A381289}. We noticed that it has  
many common terms with sequence \seqnum{A011796}. In fact, the
terms of \seqnum{A011796} coincide with the sequence $(T(n,6,2))_{n>6}$.  
Also sequence \seqnum{A011795} coincides with $(T(n,5,0))_{n>5}$.  
Both sequences, \seqnum{A011795} and \seqnum{A011796}, were created
by the first named author \cite{Broadhurst}. They arise from problems 
on multiple polylogarithms \cite{Bluemlein,Borwein}
from Feynman diagrams in quantum field theory \cite{Broadhurst2,BK} 
and this is how we began our joint work on $T(n,k,s)$, thanks to the OEIS.  

Our work is an analog, on the ring $\mathbb{Z}/n\mathbb{Z}$, of enumerating
the partitions of an integer. We found the following related, but different
kind of work. 
Maze \cite{Maze} studied the number $R(n,s,\mathcal{C})$    
of sums of the form $x_{1}+\dots+x_{k}=s$ such that the number of
repetitions of a given $x\in\mathbb{Z}/n\mathbb{Z}$ among the terms  
$x_{1},\dots,x_{k}$ is in the finite fixed set $\mathcal{C}$ of integers.
For example, the sums which are the closest to the sums we study here  
are the sums with $0$ or $1$ repetitions, i.e., $\mathcal{C}=\{0,1\}$.
In that case the number $R(n,s,\mathcal{C})$ is, with our notations, 
$R(n,s,\mathcal{C})=\sum_{k=1}^{n}T(n,k,s)$. 

The paper is structured as follows. In Section \ref{sec:A-sequence-of},  
we describe the sequence of matrices $M(k)$ and state its main properties,
some of which (e.g. its determinant and inverse) were previously known to Balandraud \cite{Balandraud}.  
 In Section
\ref{sec:T(n,k,s)-def-properties}, we study properties  
of the integers $T(n,k,s)$ and derive generating functions $G_{k,s}=\sum_{n}T(n,k,s)x^{n}$. 
Section \ref{sec:Route}
recounts the empirical route by which we arrived at the matrix construction in Theorem \ref{thm:MainIntro}. 
A proof, from representation theory, was then sent to us by Deligne \cite{Deligne}. 
It turns out that another proof due to Ramanathan was known \cite{Ramanathan}.
Section  \ref{sec:Further-remarks-and} contains further remarks and open questions.
In particular, our sequences are related to the enumeration of necklaces \cite{Hadjicostas} 
and Lyndon words \cite{Deligne,Moree}. For example, we show that $T(n,k,2)$ for $k=2\pmod4$ 
enumerates Lyndon words in an alphabet
with two letters. Section \ref{sec:Code} encodes 
the result for $T(n,k,s)$ and Section \ref{sec:OEIS} tabulates sequence numbers in
the OEIS that relate to our paper.  
Appendix A, provided by Deligne, contains his proof of Theorem \ref{thm:MainIntro}  
and suggests generalization to wider problems, dealing with virtual representations. 

\section{\label{sec:A-sequence-of}A sequence of matrices with remarkable
properties}

\subsection{Definitions of the matrices $M(k)$}

\subsubsection{Matrices indexed by finite sets}

In order to fix the notation, let us recall some well-known facts.
Let $R$ be a commutative ring with unity.
\begin{definition}
Let $I,J$ be two totally ordered finite sets. An $I\times J$-\textit{indexed 
matrix} is a function $M\in R^{I\times J}$, $M:I\times J\to R$.
Defining $M_{i,j}$ as the image of $(i,j)\in I\times J$, one may
also write $M=(M_{i,j})_{i\in I,\,j\in J}$. 
\end{definition}

If $I=J$, one speaks of an $I$-indexed matrix.  
The total order on $I$ and $J$ enables one to visualize $M$ as an array  
of $|I|$ rows and $|J|$ columns, and to compute minors. Note  that the
data of a total order on a set $I$ is equivalent to a bijection $I\to\{1,\dots,|I]|\}$.
For a set $I$, a row is a matrix indexed by $\{0\}\times I$, and
a column or vector is indexed by $I\times\{0\}$. 

If $I,J,K$ are three totally ordered finite sets and $M=(M_{i,j})_{i\in I,\,j\in J},\,N=(N_{j,k})_{j\in J,\,k\in K}$ 
are matrices indexed by $I \times J$ and $J \times K$ respectively, their product 
is the $I\times K$-indexed matrix 
\[
(\sum_{j\in J}M_{i,j}N_{j,k})_{i\in I,\,k\in K}.
\]
In particular, an $I\times J$-indexed matrix defines a linear map of  
$R$-modules $R^{I}\to R^{J}$. 

For finite sets $I,J,I',J'$, an $I\times J$ indexed matrix $M$ and  
an $I'\times J'$ indexed matrix $M'$, the Kronecker product $M\otimes M'$  
is the function $(I\times I')\times(J\times J')\to R$ defined by
\[
(M_{i,j}M'_{i',j'})_{(i,i')\in I\times I',\,(j,j')\in J\times J'}.
\]
The function $M\otimes M'$ defines a $R$-linear map $R^{I}\otimes R^{I'}\to R^{J}\otimes R^{J'}$.
Note that in general there is no canonical choice of an order on $I\times I'$
and $J\times J'$. The choice of such a total order defines a matrix.
Another choice gives another matrix and the two matrices are related
by conjugation by permutation matrices. If $M,M'$ are indexed  
by ranges of integers, there is a convention on the order of rows
and columns of $M\otimes M'$ and in general $M\otimes M'\neq M'\otimes M$
as matrices. 

\subsubsection{Divisor indexed matrices}

For an integer $k>0$, let 
\[
D_{k}=\{d\in \mathbb{N} \,|\,d\mid k\} 
\]
be its set of positive divisors, with the order induced by the usual total  
order for positive integers. We will consider $D_{k}\times D_{k}$-indexed matrices  
$M=(M_{d,d'})_{d,d'\in D_{k}},$ and $D_{k}$-indexed vectors $v=(v_{d})_{d\in D_{k}}$  
(with entries in $\mathbb{Z}$). For short, we will say that $M$  
and $v$ are $k$-indexed. For example, for $p$ prime, the  
$p$-indexed matrix 
\[
M=\left(\begin{array}{cc}
1 & 2\\
3 & 4
\end{array}\right)
\]
is such that $M_{1,1}=1$, $M_{1,p}=2$, $M_{p,1}=3$, $M_{p,p}=4$. 

If $k,k'$ are two coprime integers, then the divisors of $kk'$ are
\[
D_{kk'}=\{dd'\,|\,d\in D_{k},\,d'\in D_{k'}\}.
\]
Let $M=(M_{d,d'})_{d,d'\in D_{k}}$ and $N=(N_{d,d'})_{d,d'\in D_{k'}}$
be respectively $k$ and $k'$ indexed matrices. 
\begin{remark}
Suppose that $k,k'$ are coprime integers. The Kronecker product $M\otimes N$
of $M$ and $N$ defined by
\[
(M\otimes N)_{d_{1}d_{2},d'_{1}d'_{2}}=M_{d_{1},d_{1}'}N_{d_{2},d_{2}'}\,\,\text{for }d_{1},d_{1}'\in D_{k},\,\,d_{2},d'_{2}\in D_{k'}
\]
is a $kk'$-indexed matrix and 
\[
M\otimes N=N\otimes M.
\]
This follows from the fact that since $k$ and $k'$ are coprime,  
every divisor $d$ of $kk'$ is the product of a unique pair $(d_{1},d_{2})\in D_{k}\times D_{k'}$.  
\end{remark}

\subsection{The matrices $M(k)$}

For a prime number $p$ and a positive integer $m$, let us define the  
following $p^{m}$-indexed matrix 
\[
M(p^{m})=\left(\begin{array}{cccccc}
1 & -1 & 0 & \cdots & \cdots & 0\\
1 & e & -p & 0 & \cdots & 0\\
\vdots & e & ep & -p^{2} & \ddots & \vdots\\
\vdots & \vdots & \vdots & ep^{2} & \ddots & 0\\
\vdots & \vdots & \vdots & \vdots & \ddots & -p^{m-1}\\
1 & e & ep & ep^{2} & \cdots & ep^{m-1}
\end{array}\right),
\]
where $e=\varphi(p)=p-1$. For example 
\[
M(p)=\left(\begin{array}{cc}
1 & -1\\
1 & e
\end{array}\right),\quad M(p^{2})=\left(\begin{array}{ccc} 
1 & -1 & 0\\
1 & e & -p\\
1 & e & ep
\end{array}\right).
\]

\begin{definition}
If $k=1$, then define the matrix $M(1)=(1)$. Let $k$ be an integer
$>1$: it has a unique factorization 
\[
k=p_{1}^{m_{1}}\dots p_{j}^{m_{j}},
\]
where the $p_{1},\dots,p_{j}$ are distinct prime numbers and $m_{j}\geq1$
are integers. Let us define the matrix $M(k)$ by 
\[
M(k)=M(p_{1}^{m_{1}})\otimes\ldots\otimes M(p_{j}^{m_{j}}).
\]
\end{definition}

The matrix $M(k)$ is well defined since the Kronecker product on
divisor-indexed matrices is commutative.

\subsection{Main properties of the matrices $M(k)$}

For an integer $k>0$, let $W=W(k)$ be the $k$-indexed matrix such 
that 
\[
W_{d,d'}=\delta_{k/d,d'}\text{ for \ensuremath{d\mid k} and \ensuremath{d'\mid k},} 
\]
where
\[
\delta_{u,v}=\begin{cases}
1, & \text{if \ensuremath{u=v}};\\
0, & \text{otherwise} 
\end{cases}
\]
is the Kronecker symbol. If $A=(A_{d,d'})$ is a $k$-indexed matrix,
$WA$ is the matrix $WA=(A_{k/d,d'})_{d,d'\in D_{k}}$. 

The M\"obius function on positive integers is defined by $\mu(1)=1$  and for $n>1$ by 
\[
\mu(n)=\begin{cases}
(-1)^{m} & \text{if }n\text{ is a square-free integer with \ensuremath{m} prime factors;}\\
0 & \text{otherwise.} 
\end{cases}
\]
The Euler totient function is $\varphi(n)=n\prod_{p\mid n}(1-\frac{1}{p})$,
with a product over primes dividing $n$. 
For two integers $d,d'$ and $f$ a power of a prime number, let us
define the integer 
\[
a(d,d',f)=\begin{cases}
\varphi(\gcd(f,d')), & \text{if }\gcd(f,d')\mid\gcd(f,d);\\
\mu\left(\frac{\gcd(f,d')}{\gcd(f,d)}\right)\gcd(f,d), & \text{otherwise.} 
\end{cases}
\]
The main properties of the matrices $M(k)$ are summarized by the following theorem. 
\begin{theorem}
\label{thm:Merveilles}For all positive integers $k$, the following 7 properties of $M(k)$ hold. 
\begin{enumerate}
\item For an integer $k$, let $F_{k}$ be the set of integers of the form
$p^{e}$ where $p$ is a prime dividing $k$ and $e>0$ is the exponent 
of $p$ in $k$. For two divisors $d,d'$ of $k$, the $d,d'$ entry
of $M(k)$ is 
\begin{equation}
M(k)_{d,d'}=\prod_{f\in F_{k}}a(d,d',f).\label{eq:BelleFormule}
\end{equation}
\item $\left(W(k)M(k)\right)^{2}=kI$, where $I$ is the identity matrix.
\item $\det(M(k))=k^{|D_{k}|/2}$. 
\item Consider $w=(1)_{d\mid k}$. The vector $M(k)w$ (the sum of the columns)
is equal to the vector $v=(k\delta_{d,k})_{d\mid k}$. 
\item The first row of $M(k)$ is $M_{1,d}=\mu(d)$, $\forall d\mid k$. 
\item The last row of $M(k)$ is $M_{k,d}=\varphi(d)$, $\forall d\mid k$. 
\item The first column of $M(k)$ is $M_{d,1}=1$, $\forall d\mid k$. 
\end{enumerate}
\end{theorem}

\begin{proof}
Let us prove Theorem \ref{thm:Merveilles}.
\begin{enumerate}
\item It is easy to check that Formula \eqref{eq:BelleFormule} is true 
for $k=p^{e}$, a power of a prime. Suppose that $k=k_{1}k_{2}$ is 
a product of two coprime integers. By induction the result is true
for $k_{1}$ and $k_{2}$: 
\[
M(k_{1})_{d_{1},d'_{1}}=\prod_{f\in F_{k_{1}}}a(d_{1},d'_{1},f),\,\,M(k_{2})_{d_{2},d'_{2}}=\prod_{f\in F_{k_{2}}}a(d_{2},d'_{2},f).
\]
{}From the construction of $M(k)$, one has 
\[
M(k)_{d_{1}d_{2},d_{1}'d_{2}'}=M(k_{1})_{d_{1},d'_{1}}M(k_{2})_{d_{2},d'_{2}}.
\]
The set $F_{k}$ is the disjoint union of $F_{k_{1}}$ and $F_{k_{2}}$.
Suppose that $f\in F_{k_{1}}$. Then $\gcd(d_{1}d_{2},f)=\gcd(d_{1},f)$ 
and $\gcd(d_{1}'d_{2}',f)=\gcd(d_{1}',f)$, and from the definition
of $a$, one has
\[
a(d_{1}d_{2},d_{1}'d_{2}',f)=a(d_{1},d_{1}',f),
\]
and the analogous result holds when $f\in F_{k_{2}}.$ This proves that  
\[
M(k)_{d_{1}d_{2},d_{1}'d_{2}'}=\prod_{f\in F_{k}}a(d_{1}d_{2},d_{1}'d_{2}',f).
\]
\item When $k$ is a power of a prime, a direct computation gives  $\left(W(k)M(k)\right)^{2}=kI$.  
For the general case, $k=p_{1}^{e_{1}}\dots p_{m}^{e_{m}}$, we remark
that $W(k)=\otimes_{j=1}^{m}W(p_{j}^{e_{j}})$ and  
\[
W(k)M(k)=\otimes_{j=1}^{m}W(p_{j}^{e_{j}})M(p_{j}^{e_{j}}). 
\]
Moreover for $k$-indexed (respectively $k'$-indexed) matrices $A,A'$
(respectively $B,B'$) with coprime integers $k,k'$, the multiplication
(Hadamard product) is such that: $(A\otimes B)(A'\otimes B')=AA'\otimes BB'$.
This implies that $\left(W(k)M(k)\right)^{2}=kI$ for every $k$.  
\item When $k=p^{e}$ is a power of a prime, a direct computation by subtracting  
rows yields $\det(M(p^{e}))=p^{e(e+1)/2}$. We then use the property 
that for square matrices $A,B$ of size $m,n$ respectively, one has
$\det(A\otimes B)=\det(A)^{n}\det(B)^{m}$. 
\item The claim follows from direct computation on $M(p^{e})$ for a prime $p$ 
and $e\geq1$. Suppose that the result is known for integers $<k$
and that $k=k_{1}k_{2}$ is a product of two coprime integers. Then
\begin{align*} 
\sum_{d_{1}\mid k_{1},\,d_{2}\mid k_{2}}M(k)_{d_{1}d_{1}',d_{2}d_{2}'} & =\sum_{d_{2}\mid k_{2}}\sum_{d_{1}\mid k_{1}}M(k_{1})_{d_{1},d_{1}'}M(k_{2})_{d_{2},d_{2}'}
=\sum_{d_{2}\mid k_{2}}k_{1}\delta_{k_{1},d_{1}'}M(k_{2})_{d_{2},d_{2}'}\\
& =k_{1}\delta_{k_{1},d_{1}'}k_{2}\delta_{k_{2},d_{2}'}=k\delta_{k,d_{1}'d_{2}'}. 
\end{align*}
\end{enumerate}
Points (5), (6) and (7) are direct consequences of Formula \eqref{eq:BelleFormule}.
\end{proof}
We also note the following result when $k$ is not divisible by the square of a prime.
\begin{proposition}
If $k$ is square-free, then  
\[
M(k)_{d,d'}=\mu\left(\frac{d'}{\gcd(d,d')}\right)\varphi(\gcd(d,d')).
\]
\end{proposition}
\begin{proof}
We are in the case where every exponent of the factorization of $k$ is unity, i.e., $F_k$ is the set of prime numbers dividing $k$. 
Let $p\in F_k$ such that $p$ divides $d'/\gcd(d,d')$. Then $\gcd(p,d')=p$ does not divide 
$\gcd(p,d)=1$, and  by definition, we have $a(d,d',p)=\mu(\gcd(p,d'))$. 
Let $p\in F_k$ such that $p$ divides $\gcd(d,d')$. Then  
$a(d,d',p)=\varphi(\gcd(p,d'))$.
For the remaining primes $p\in F_k$, we have $a(d,d',p)=1$. Therefore, using the
multiplicativity of the functions $\mu$ and $\varphi$, and applying Formula \eqref{eq:BelleFormule}, we obtain
$M(k)_{d,d'}=\prod_{p|k}a(d,d',p)=\mu\left(\frac{d'}{\gcd(d,d')}\right)\varphi(\gcd(d,d'))$.
\end{proof}

\section{\label{sec:T(n,k,s)-def-properties}The sequences $(T(n,k,s))_{n\geq k}$} 

\subsection{Definitions and properties of the integers $T(n,k,s)$}  

Let us recall that for integers $n,k,s$, with $n\geq k>0$, we denote by
$S(n,k,s)$ the set of order $k$ subsets $Q$ of $\mathbb{Z}/n\mathbb{Z}$
such that the sum $\sum_{a\in Q}a$ equals $s$ modulo $n$. Moreover,
the order of $S(n,k,s)$ is denoted by 
\[
T(n,k,s)=|S(n,k,s)|.
\]
We now  prove the following result. 
\begin{proposition}
\label{prop:action}The function $s\to T(n,k,s)$ depends only on
the class of $s$ modulo $\gcd(n,k)$. Let $\beta\in(\mathbb{Z}/\gcd(n,k)\mathbb{Z})^{*}$.
Then $T(n,k,s)=T(n,k,\beta s)$ for all $s\in\mathbb{Z}/\gcd(n,k)\mathbb{Z}$. 
\end{proposition}

\begin{proof}
One has the relation $S(n,k,s)=S(n,k,s+n)$, therefore one may consider
that $s$ is in $\mathbb{Z}/n\mathbb{Z}$. Moreover for every $Q=\{a_{1},\dots,a_{k}\}\in S(n,k,s)$
and every $ \alpha\in\mathbb{Z}/n\mathbb{Z}$, the set 
\[
Q( \alpha):=\{a_{1}+ \alpha,\dots,a_{k}+ \alpha\}
\]
 is in $S(n,k,s+k \alpha)$. Also we remark that if $\beta\in(\mathbb{Z}/n\mathbb{Z})^{*}$
is invertible, then 
\[
\beta Q:=\{\beta a_{1},\dots,\beta a_{k}\}
\]
is in $S(n,k,\beta s)$.  For $\alpha\in\mathbb{Z}/n\mathbb{Z}$ and $\beta\in(\mathbb{Z}/n\mathbb{Z})^{*}$, the maps $Q\to Q( \alpha)$
and $Q\to\beta Q$ 
are bijections, which implies that 
\begin{align*} 
T(n,k,s)&=T(n,k,s+k\alpha),\;\forall\alpha\in\mathbb{Z}/n\mathbb{Z},\\ 
T(n,k,s)&=T(n,k,\beta s),\;\forall\beta\in(\mathbb{Z}/n\mathbb{Z})^{*}. 
\end{align*} 
Since also $T(n,k,s)=T(n,k,s+n)$, the function 
\[
s\to T(n,k,s)
\]
is periodic of period $m=\gcd(n,k)$.  
Thus we may consider $s$ in $\mathbb{Z}/m\mathbb{Z}$.  
Then the action of $(\mathbb{Z}/n\mathbb{Z})^{*}$ factors through
the surjective map
\[
(\mathbb{Z}/n\mathbb{Z})^{*}\to(\mathbb{Z}/m\mathbb{Z})^{*},  
\]
by which we mean that the following relation holds:
\[
T(n,k,s)=T(n,k,\beta s),\;\forall\beta\in(\mathbb{Z}/m\mathbb{Z})^{*}. 
\]
\end{proof}
\begin{remark}\label{Remaction}
The proof of Proposition \ref{prop:action} shows that the affine
group $( \mathbb{Z}/n \mathbb{Z})^{*}\rtimes \mathbb{Z}/n \mathbb{Z}$ acts on the set $\{S(n,k,s)\,|\,s\in \mathbb{Z}/n \mathbb{Z}\}$. 
\end{remark}

\begin{corollary}\label{cor:complement} 
(a) For $n\geq k>0$, 
\[T(n,k,s)=T(n,k,k+s)=T(n,k,k-s).\] 
(b) For $n> k>0$, 
\[T(n,k,s)=\begin{cases}T(n,n-k,s),&\text{if }n\text{ is odd};\\ 
T(n,n-k,s+n/2),&\text{if }n\text{ is even}.\end{cases}\] 
\end{corollary}

\begin{proof}
The first part is immediate from Proposition \ref{prop:action}. For
the second part, we remark that the complement of a subset of order 
$k$ has order $n-k$, and we use the fact that the sum of elements
of $\mathbb{Z}/n\mathbb{Z}$ is $0$ if $n$ is odd and $n/2$ if
$n$ is even.
\end{proof}

\subsection{Generating functions}

We know that for integers $s,t$ such that $\gcd(k,s)=\gcd(k,t)$,
one has $T(n,k,s)=T(n,k,t)$ for all integers $n$. Therefore, for
a fixed integer $k$, to compute the value $T(n,k,s)$, it suffices
to consider $s$ in $D=D_{k}$, the ordered set of positive divisors of $k$.  
Let $P=P_{k}$ be the $k$-indexed vector $P_{k}=(P_{k,d})_{d\in D}$,
where for $d\in D$, $P_{k,d}$ is the formal power series
\[
P_{k,d}=\frac{(-1)^{k-k/d}}{k}\frac{x^{k}}{(1-x^{d})^{k/d}}.
\]
For $d\mid k$, let us define the formal power series 
\[
G_{k,d}=\sum_{n\geq k}T(n,k,d)x^{n}, 
\]
and let $G_{k}$ be the vector $G_{k}=(G_{k,d})_{d\in D}$. A direct
computation shows that Theorem \ref{thm:MainIntro} can be rephrased
as follows. 
\begin{theorem}
\label{thm:MAIN}The matrix $M(k)=(M_{d,d'})_{d,d'\in D}$ is such
that
\[
G_{k}=M(k)P_{k}.
\]
\end{theorem}

\section{\label{sec:Route}Route to proof}  

\subsection{Last row of $M(k)$ and generating function of $T(n,k,0)$}
In a comment on \seqnum{A032801}, Hadjicostas remarked that Barnes proved the following result. 
\begin{proposition}
\label{prop:(Barnes-)}(Barnes \cite[Lemma 5.1]{Barnes}). The number of elements of $S(n,k,0)$ is  
\[
T(n,k,0)=\frac{1}{n}\sum_{s\mid\gcd(n,k)}(-1)^{k-k/s}\varphi(s)\binom{n/s}{k/s}. 
\]
\end{proposition}

We recall that Theorem \ref{thm:MAIN} states the equality $G_{k}=M(k)P_{k}$.  
The following result is a consequence of Barnes's formula in Proposition \ref{prop:(Barnes-)}. 

\begin{proposition}
\label{prop:DB2}Theorem \ref{thm:MAIN} is satisfied for the
index $k$ of vector $G_{k}$. In other words:
\[
G_{k,k}=\sum_{s\mid k}\frac{(-1)^{k-k/s}\varphi(s)}{k}\frac{x^{k}}{(1-x^{s})^{k/s}}=\left(M(k)P_{k}\right)_{k}.
\]
\end{proposition}

\begin{proof}  
From Barnes's formula in Proposition \ref{prop:(Barnes-)}, one has
\[
T(n,k,k)=T(n,k,0)=\frac{1}{n}\sum_{s\mid\gcd(n,k)}(-1)^{k-k/s}\varphi(s)\binom{n/s}{k/s}. 
\]
For $d\mid k$, let $ \mathbb{N}_{d,k}=\{n\in \mathbb{N}\,|\,\gcd(n,k)=d\}$.  Then 
\[
G_{k,k}=\sum_{n\geq k}T(n,k,k)x^{n}=\sum_{d\mid k}\sum_{s\mid d}\sum_{n\in \mathbb{N}_{d,k}}
(-1)^{k-k/s}\varphi(s)\binom{n/s}{k/s}\frac{x^n}{n}.  
\]
For $s\mid k$, one has $(\cup_{s\mid d,\,d\mid k} \mathbb{N}_{d,k})=s \mathbb{N}$ and hence 
\[
G_{k,k}=\sum_{s\mid k}\frac{(-1)^{k-k/s}\varphi(s)}s\sum_{m\geq1}\binom{m}{k/s}\frac{x^{sm}}{m}  
=\sum_{s\mid k}\frac{(-1)^{k-k/s}\varphi(s)}{k}\frac{x^{k}}{(1-x^{s})^{k/s}} 
\]
where the last equality is obtained from the formal power series 
\[
\frac{y^j}{(1-y)^j}=j\sum_{m\geq 1}\binom{m}{j}\frac{y^m}{m} 
\]
with $j=k/s\in\mathbb{N}$, and $y=x^s$. 
\end{proof} 

\subsection{A route to the sequence of matrices $M(k)$}  

Here, we recount how we arrived at Theorem \ref{thm:MainIntro}, 
and its equivalent formulation by generating functions in Theorem \ref{thm:MAIN}. 
The generating functions of the sequences \seqnum{A011795} and \seqnum{A011796} 
(which are related to $T(n,5,1)$ and $T(n,6,2)$) are known. With
Proposition \ref{prop:DB2}, this makes plausible the following 
less specific proposition.
\begin{proposition}
\label{thm:Light}The generating functions $G_{k,d},\,d\mid k$ are
integral combinations of 
\[
P_{k,d'}=\frac{(-1)^{k-k/d'}}{k}\frac{x^{k}}{(1-x^{d'})^{k/d'}} 
\]
with $G_{k,d}=\sum_{d'\mid k}M_{d,d'}P_{k,d'}$ and coefficients $M_{d,d'}\in \mathbb{Z}$. 
\end{proposition}

To exhibit how this led us to the more precise Theorem \ref{thm:MAIN}, we 
prove the following proposition. 
\begin{proposition}
\label{prop:DB1}Suppose that Proposition \ref{thm:Light} holds. Then 
\[
\forall d\mid k,\,\begin{cases}
\sum_{d'\mid k}M_{d,d'}=k\delta_{k,d}, & \text{if \ensuremath{k} is odd;}\\ 
\sum_{d'\mid k}(-1)^{k/d'}M_{d,d'}=k\delta_{k/2,d}, & \text{if \ensuremath{k} is even.} 
\end{cases}
\]
\end{proposition}

\begin{proof}
Let $d\mid k$ be a divisor of $k$. By definition, $\frac{1}{x^{k}}G_{k,d}=\sum_{n\geq k}T(n,k,d)x^{n-k}$, and 
thus $\lim_{x\to0}\frac{1}{x^{k}}G_{k,d}=T(k,k,d)$. Therefore, if 
$G_{k,d}=\sum_{d'\mid k}M_{d,d'}P_{k,d'}$, 
we obtain 
\[
T(k,k,d)=\frac1k\sum_{d'\mid k}(-1)^{k-k/d'}M_{d,d'}. 
\]
For odd $k$, $T(k,k,d)=0$ if $d$ is not $0\bmod k$, 
while for even $k$, $T(k,k,d)=0$ if $d$ is not $k/2\bmod k$. Moreover, 
$T(k,k,k)=1$ if $k$ is odd, and $T(k,k,k/2)=1$ if $k$ is even.  
\end{proof}
Proposition \ref{prop:DB1} determines the final column of $M=M(k)$
and Proposition \ref{prop:DB2} determines its final row. Hence the 
matrix $M(k)$ is determined by its remaining $(|D_k|-1)^{2}$ elements. 
These elements may be computed from values of $T(n,k',s')$ with $n/2\leq k'<k\leq n$,  
using Corollary \ref{cor:complement}, as follows. 
By definition, 
\[
G_{k,s}=\sum_{n=k}^{2k-1}T(n,k,s)x^n+O(x^{2k}). 
\]
For $n\in\{k,\dots,2k-1\}$, one has $n-k\in\{0,\dots,k-1\}$. Moreover, 
\[
T(n,k,s)=T(n,n-k,s+\epsilon_{n}),
\]
where $\epsilon_{n}=n/2$ if $n$ is even, and $\epsilon_{n}=0$ if $n$ is odd.
Thus knowledge of the values $T(n,k',s')$ for $k'<k$, $s'\in\mathbb{Z}/n\mathbb{Z}$ 
enables one to compute the first $k$ terms of the series $G_{k,s}$. 
If $G_{k,s}$ is as claimed in Proposition \ref{thm:Light}, i.e.,
an integral combination of the series $P_{k,d},\,d\mid k$, then the 
first $k$ coefficients of $G_{k,s}$ determine the desired 
integers $M_{s,d}$. Indeed, the equality 
\begin{equation}
\sum_{n=k}^{2k-1}T(n,k,s)x^{n}=\sum_{d\mid k}M_{\gcd(k,s),d}P_{k,d}+O(x^{2k})\label{eq:experiment} 
\end{equation}
gives a linear system, with unknowns $M_{t,d}$ for $t\mid k,\,d\mid k$, 
satisfying $k^{2}$ equations.  

The first $10000$ matrices of the sequence $(M(k))_{k\geq1}$ have
been obtained that way: by supposing that Proposition \ref{thm:Light}
holds true, and by resolving Equation \eqref{eq:experiment}. Note 
that there are more equations than unknowns. Yet each time we found 
a unique solution. 

After determining these $10000$ matrices, a pattern appeared for $M(p^{e})$, for  
$p$ a prime number, and we discovered that one may obtain the matrix
$M(k)$ by factorizing $k$ and taking the Kronecker product over
that factorization.
Thus our work subsumed a vast amount of empirical observation. 

\subsection{Proof of the matrix construction using Ramanujan sums}

In an initial version of this paper, we presented Theorem \ref{thm:MainIntro}
as a conjecture, supported by a large amount of testing and some partial proofs.
We soon received helpful messages, from  E.\ Balandraud, P.\ Deligne, and P.\ Hadjicostas, 
which enable one to prove Theorem \ref{thm:MainIntro}, as follows.

\begin{proposition}\label{prop:rs}
Formulas \eqref{Muv} and \eqref{Mtd} define the Ramanujan sum
\[ 
M(k)_{s,d} = c_d(s) := \sum_{j\in(\mathbb{Z}/d\mathbb{Z})^{*}}\zeta_d^{js}
=\mu\left(\frac{d}{\gcd(s,d)}\right)\frac{\varphi(d)}{\varphi\left(\frac{d}{\gcd(s,d)}\right)}
\]
where $s,d$ are positive divisors of $k>0$ and $\zeta_d=\exp(2\pi i/d)$ is a primitive $d$-th root of unity.
\end{proposition}

\begin{proof}
The equality between Ramanujan's sum $c_d(s)$ and 
$\mu\left(\frac{d}{\gcd(s,d)}\right)\frac{\varphi(d)}
{\varphi\left(\frac{d}{\gcd(s,d)}\right)}$ 
is proved, for example, by Hardy and Wright \cite[Theorem 272, p.\ 238]{HW}. 
Formula \eqref{Muv} establishes that $M(k)_{s,d}=c_d(s)$ 
when $k$ is a prime power. The general result then follows 
from the multiplicativity \cite[Theorem 67, p.\ 56]{HW} of Ramanujan  
sums, with $c_{dd'}(s)=c_d(s)c_{d'}(s)$ when $d$ and $d'$ 
are coprime. If $k$ and $k'$ are coprime, then 
$\gcd(d,d')=1$ when $d\mid k$ and $d'\mid k'$.
\end{proof}

We were alerted to Proposition \ref{prop:rs} by Balandraud, who has studied Ramanujan
sums and some properties of the matrix $M(k)$ \cite[Proposition 2]{Balandraud},
and by Hadjicostas, who has added important notes to the OEIS entries \seqnum{A054535} and \seqnum{A267632}.
Proposition \ref{prop:rs} shows that Theorem \ref{thm:MainIntro} gives a matrix construction that is
equivalent to a theorem by Ramanathan \cite[Theorem 4]{Ramanathan}, 
proved 80 years ago.  Moreover, Hadjicostas suggested the following novel 
proposition \cite{Hadjicostas2}, whose proof we now provide. 

\begin{proposition}
 For $n\geq k>0$ and $s\mid k$, one has 
\begin{equation}T(n,k,s)=\sum_{d\mid\gcd(n,s)}(-1)^{k-k/d}T(n/d,k/d,1).
\label{PH}\end{equation}
\end{proposition}

\begin{proof} 
Abbreviating $m:=\gcd(n,k)$  and applying Theorem \ref{thm:MainIntro} to $T(n/d,k/d,1)$ 
on the  r.h.s.\ of Formula \eqref{PH}, one obtains
\[R(n,k,s):=\sum_{d\mid\gcd(n,s)}(-1)^{k-k/d}T(n/d,k/d,1)=
\sum_{d\mid\gcd(n,s)}\sum_{t|(m/d)}\frac{d}{n}(-1)^{k-k/(td)}\mu(t)\binom{n/(td)}{k/(td)}.\]
With $t=e/d$, this gives
\[R(n,k,s)=\frac{1}{n}\sum_{e|m}(-1)^{k-k/e}\binom{n/e}{k/e}
\sum_{d\mid\gcd(s,e)}\mu\left(\frac{e}{d}\right)d.\]
The sum over $d$, with $d\mid s\mid k$,
and $d\mid e\mid m \mid n$, 
gives the Ramanujan sum $c_e(s)$, as recorded by Hardy and Wright \cite[Theorem 271, p.\ 237]{HW}.
Hence we obtain
\[R(n,k,s)=\frac{1}{n}\sum_{e|m}(-1)^{k-k/e}\binom{n/e}{k/e}M(k)_{s,e}
=T(n,k,s)\]
from Proposition \ref{prop:rs}, followed by Theorem \ref{thm:MainIntro}.
\end{proof} 
 
For Deligne's proof of our matrix construction and his further remarks, see Appendix A,
where Conjecture 1 refers to what is now established in Theorem \ref{thm:MainIntro}.

\section{\label{sec:Further-remarks-and}Further remarks and open questions}

\subsection{Necklaces and sequences $T(n,k,0),$ $T(n,k,k/2)$}

A necklace on $n$ beads, with $k$ black beads and $n-k$ white
beads, can be seen as an orbit $O(Q)$, of a subset $Q\subset\mathbb{Z}/n\mathbb{Z}$ 
with $k$ elements, under the action of $\mathbb{Z}/n\mathbb{Z}$ defined by  
\[
r+Q:=\{r+q\,|\,q\in Q\},\,\text{for }r\in\mathbb{Z}/n\mathbb{Z}.
\]
Burnside's theorem \cite[Chapter 6]{Steinberg} gives the number $N(n,k)$ of necklaces, i.e., the 
number of orbits, as 
\[
N(n,k)=\frac{1}{n}\sum_{d\mid\gcd(n,k)}\varphi(d)\binom{n/d}{k/d}.
\]

The \emph{last} row of $M(k)$ is given $M(k)_{k,d}=\varphi(d)$, for $d\mid k$.   
Thus for odd $k$, one has 
\[
T(n,k,0)=T(n,k,k)=\frac{1}{n}\sum_{d\mid\gcd(n,k)}(-1)^{k-k/d}\varphi(d)\binom{n/d}{k/d}=N(n,k) 
\]
since, for odd $k$, one has $(-1)^{k-k/d}=1$, for all $d\mid k$. 

When $k$ is even, Formulas \eqref{Muv} and \eqref{Mtd} give the \emph{penultimate} row of $M(k)$ as 
\[
M(k)_{k/2,d}=(-1)^{k/d}M(k)_{k,d}=(-1)^{k/d}\varphi(d), 
\]
for all $d\mid k$. Hence Theorem \ref{thm:MainIntro} gives
\[
T(n,k,k/2)=\frac{1}{n}\sum_{d\mid\gcd(n,k) }(-1)^{k-k/d}M(k)_{k/2,d}\binom{n/d}{k/d}
=\frac{1}{n}\sum_{d\mid\gcd(n,k)}\varphi(d)\binom{n/d}{k/d}=N(n,k). 
\]
We thus prove the following result.
\begin{proposition}\label{pr:necklace}
For $n>k>0$, the number 
of necklaces with $k$ black beads and $n-k$ white
beads is $N(n,k)=T(n,k,0)$ if $k$ is odd, and $N(n,k)=T(n,k,k/2)$
if $k$ is even. 
\end{proposition}

Thus we see that these two enumerative problems:  
\begin{description}
\item [{a}] counting the number of necklaces, i.e., the number of orbits
$O(Q)$ of order $k$ subsets $Q$ of $\mathbb{Z}/n\mathbb{Z}$ under 
translations, and 
\item [{b}] counting the number of order $k$ subsets $Q$ of $\mathbb{Z}/n\mathbb{Z}$ 
with sum equal to $0$ or $k/2$ (according as whether $k$ is odd or even)
\end{description}
share the same formula.

Consider the map $\gamma$ that assigns, to each order $k$ subset $Q$ 
of $\mathbb{Z}/n\mathbb{Z}$, its orbit $O(Q)$. 
For $k$ odd, respectively even, it is natural to ask whether $S(n,k,k)$, 
respectively $S(n,k,k/2)$, is sent by $\gamma$ bijectively to the set
of necklaces. If so, this would immediately explain the coincidence of these two enumerations. 

For example, suppose that $k$ and $n$ are coprime. The translation 
by $r\in\mathbb{Z}/n\mathbb{Z}$ maps $S(n,k,s)$ to $S(n,k,s+kr)$.
Since $k$ is coprime to $n$, that action is faithful on the sets
$S(n,k,s),\,s\in\mathbb{Z}/n\mathbb{Z}$. Thus, for fixed $s\in\mathbb{Z}/n\mathbb{Z},$
the set $S(n,k,s)$ is in bijection with the set of necklaces via $\gamma$.

On the contrary, in many examples that we tested with $\gcd(n,k)\neq1$, 
the map $\gamma$ is not a bijection, and hence we lack an immediate explanation for why 
the two formulas coincide. 

\subsection{Lyndon words and sequences $T(n,k,1)$, $T(n,k,2)$}

An \emph{aperiodic necklace} (equivalently, a \emph{binary Lyndon word}) with $k$ 
black beads and $n-k$ white beads is a necklace $O(Q)=\{r+Q\,|\,r\in\mathbb{Z}/n\mathbb{Z}\}$,
with $Q\subset\mathbb{Z}/n\mathbb{Z}$ of order $k$, such that 
the order of the orbit $O(Q)$ is $n$. It is well-known (see comments on \seqnum{A051168}) that the 
number of such binary Lyndon words is  
\[
L(n,k)=\frac{1}{n}\sum_{d\mid\gcd(n,k)}\mu(d){n/d \choose k/d}. 
\]
This formula is obtained by applying M\"obius inversion to 
${n \choose k}=\sum_{d\mid\gcd(n,k)}\frac{n}{d} L(\frac{n}{d},\frac{k}{d})$.  

Therorem \ref{thm:MainIntro} gives the \emph{first} row of matrix $M(k)$  
as $M(k)_{1,d}=\mu(d),\,d\mid k$, and hence 
\[
T(n,k,1)=\frac{1}{n}\sum_{d\mid\gcd(n,k)}(-1)^{k-k/d}\mu(d){n/d \choose k/d}. 
\]
For $k$ odd and $d\mid k$, one has $(-1)^{k-k/d}=1$. For $k=0\pmod4$ 
and $d\mid k$, one has $(-1)^{k-k/d}\mu(d)=\mu(d)$. Hence $T(n,k,1)=L(n,k)$, in both of these cases. 

In the notable remaining case, with $k=2\pmod4$, the \emph{second} row of matrix $M(k)$ is 
\[
M(k)_{2,d}=\begin{cases}\mu(d),&\text{if }d\text{ is odd;}\\ 
\mu(d/2),&\text{if }d\text{ is even,}\end{cases} 
\]
for all $d\mid k$. Hence $M(k)_{2,d}=(-1)^{k/d}\mu(d)$ and one has the relation 
\[
T(n,k,2)=\frac{1}{n}\sum_{d\mid\gcd(n,k)}(-1)^{k-k/d}M(k)_{2,d}{n/d \choose k/d} 
=\frac{1}{n}\sum_{d\mid\gcd(n,k)}\mu(d){n/d \choose k/d}=L(n,k).  
\]
Thus we have proved the following proposition. 

\begin{proposition}\label{pr:Lyndon}
For $n>k>0$, the number  of aperiodic necklaces with $k$ black beads and $n-k$ white 
beads is $L(n,k)=T(n,k,2)$ if $k=2\pmod4$, and $L(n,k)=T(n,k,1)$ otherwise. 
\end{proposition}

\begin{remark} 
As noted in the introduction, our work began by observing coincidences of terms 
in sequence \seqnum{A011796},  for $(L(n,6))_{n>6}$, with evaluations of $T(n,6,2)$. 
These are now proved by Proposition \ref{pr:Lyndon}. 
Yet we lack an immediate explanation of the established equality $L(n,6)=T(n,6,2)$. 
\end{remark} 

\subsection{On the matrices $M(p^{e})$ for $p$ a prime number}

Whenever $k=k_{1}k_{2}$ with coprime integers $k_{1},k_{2}$,
the relation $M(k)=M(k_{1})\otimes M(k_{2})$ holds. It 
is tempting to ask if there could be such results when $k=p^{e}$
is a power of a prime $p$, with the Kronecker product replaced by
a law $\star$ such that the product $M(p^{e})\star M(p^{e'})$ of 
$M(p^{e})$, of size $e+1$, and indexed by $D_{p^{e}}$, times $M(p^{e'})$, 
of size $e'+1$, is the matrix $M(p^{e+e'})$ of size $e+e'+1$ indexed
by $D_{p^{e+e'}}$. In particular, one would like to have $M(p)^{\star e}=M(p^{e})$. 
Let $f_{1},f_{p}$ be the canonical base of $R^{2}$. The symmetric
product $S^{e}M(p)$ of $M(p)$ is obtained as the action of $M(p)^{\otimes e}$
on the quotient of $(R^{2})^{\otimes e}$ by the relations 
\[
f_{j_{1}}\otimes\dots\otimes f_{j_{e}}=f_{j_{1}'}\otimes\dots\otimes f_{j_{e}'}
\]
if and only if $\prod_{k=1}^{e}j_{k}=\prod_{k=1}^{e}j_{k}'$. It is
then natural to enquire whether $M(p^{e})$ is a conjugate of the symmetric 
product $S^{e}M(p)$ of $M(p)$. Interestingly, this is not the case: 
the characteristic polynomial of 
\[
S^{2}M(p)=\left(\begin{array}{ccc}
1 & 1 & 1\\
-2 & p-2 & 2(p-1) \\
1 & 1-p & (p-1)^2
\end{array}\right)
\]
is 
\[
(X-p)(X^{2}+(2p-p^{2})X+p^{2})
\]
whereas the characteristic polynomial of $M(p^{2})$ is 
\[
(X-p)(X^{2}+(p-p^{2})X+p^{2}).
\]
In consequence,  we regard our formula for $M(p^e)$ as primitive. 
There appears to be no simple way of deriving it from the result for $M(p)$. 

\section{\label{sec:Code}Code}  

We used this compact statement of Theorem \ref{thm:MainIntro}
in {\tt Pari/GP}:
\begin{verbatim}
{T(n,k,s)=local(q=factor(k),t=gcd(k,s),u,v,w);
q=vector(#q~,i,q[i,1]^q[i,2]);
1/n*sumdiv(gcd(n,k),d,w=(-1)^(k-k/d)*binomial(n/d,k/d);
for(i=1,#q,u=gcd(q[i],t);v=gcd(q[i],d);
w*=if(v>u,moebius(v/u)*u,eulerphi(v)));w);}
\end{verbatim}
and its translation to {\tt Magma}:
\begin{verbatim}
function T(n, k, s)
q := [f[1]^f[2] : f in Factorization(k)]; 
t := GCD(k, s); r := 0;
for d in Divisors(GCD(n, k)) do
w := (-1)^(k - k div d) * Binomial(n div d, k div d);
for i in [1..#q] do
u := GCD(q[i], t); v := GCD(q[i], d);
w *:= v gt u select MoebiusMu(v div u) * u else EulerPhi(v);
end for; r +:= w; end for; return r div n; end function;
\end{verbatim}

\section{\label{sec:OEIS}OEIS sequences}

We give here pairs $(k,0)$ (respectively  $(k,s)$, $s\neq 0$) for which Proposition \ref{prop:DB2}
(respectively Theorem \ref{thm:MainIntro}) 
gives a sequence $(T(n,k,s))_{n\geq k}$ in the OEIS:\\
$(3,0)$: \seqnum{A007997}, $(3,1)$: \seqnum{A001840},\\
$(4,0)$: \seqnum{A032801}, $(4,1)$: \seqnum{A006918}, $(4,2)$: \seqnum{A008610},\\ 
$(5,0)$: \seqnum{A008646}, $(5,1)$: \seqnum{A011795},\\
$(6,0)$: \seqnum{A381289}, $(6,1)$: \seqnum{A381290}, $(6,2)$: \seqnum{A011796}, $(6,3)$: \seqnum{A032191},\\  
$(7,0)$: \seqnum{A032192}, $(7,1)$: \seqnum{A011797},\\
$(8,0)$: \seqnum{A381291}, $(8,1)$: \seqnum{A031164}, $(8,2)$: \seqnum{A381350}, $(8,4)$: \seqnum{A032193},\\  
$(9,0)$: \seqnum{A032194}, $(9,1)$: \seqnum{A263318}, $(9,3)$: \seqnum{A381351},\\
on the understanding that the index $n$ is suitably offset. 
The cases with $s=0$ are the subject of the triangular array \seqnum{A267632}. 
Binary Lyndon words are the subject of the triangular array \seqnum{A051168}. 
Ramanujan sums are the subject of the square array \seqnum{A054535}.
We obtained \seqnum{A379920}, for irreducible conics,  
by eliminating products of lines from \seqnum{A381289}. 

\section{Acknowledgments} We have benefited greatly from discussions with 
\'Eric Balandraud, Pierre Deligne, and Petros Hadjicostas. Encouragement and advice came from 
Johannes  Bl\"umlein,  Dirk Kreimer, Pieter Moree,  Michael Oakes, Tanay Wakhare, and Doron Zeilberger. 

\appendix

\renewcommand{\thetheorem}{\Alph{section}\arabic{theorem}}
\renewcommand{\thecorollary}{\Alph{section}\arabic{corollary}}
\renewcommand{\thelemma}{\Alph{section}\arabic{lemma}}
\renewcommand{\theproposition}{\Alph{section}\arabic{proposition}}

\counterwithin*{theorem}{section} 

\newtheorem*{example*}{Example}
\setcounter{theorem}{0}

\section{Appendix: Proof of Conjecture 1 by Pierre Deligne}

For $A$ a finite abelian group, $k \geq 0$ an integer and $a \in A$, define $T(A, k, a)$ to be the number 
of $k$-element subsets of $A$ summing up to $a$. We will prove a formula for $T$ which, 
for $A$ cyclic, reduces to Conjecture 1.

The Pontryagin dual $A^* := \text{Hom}(A, \mathbb{C}^*)$ of $A$ is isomorphic to $A$.  
Rather than $T(A, k, a)$, it will be convenient to consider $T(A^*, k, a^*)$, for $a^* \in A^*$. 
We will identify $A^*$ with the set of isomorphism classes of one-dimensional, equivalently irreducible, 
complex representations of $A$.

Let $R$ be the regular representation of $A$: the space of (complex) functions on $A$, on 
which $A$ acts by translations. We will use two bases of $R$: the function $\delta_a$: $1$ 
at $a$, $0$ elsewhere ($a \in A$), and the functions $a^*$ ($a^* \in A^*$). An element $a$ 
of $A$ acts on the function $a^*$ by multiplying it by $a^*(a)^{-1}$.

\begin{proposition} \label{propositionA1}
The number $T(A^*, k, a^*)$ is the multiplicity with which the representation $a^*$ occurs in $\wedge^k R$.
\end{proposition}
  
\begin{proof}
As $R$ is the direct sum of the representations $b^*$ ($b^* \in A^*$), $\wedge^k R$ is the direct sum, 
over the $k$-element subsets $B$ of $A^*$, of the representations $\otimes_B b^*$.
\end{proof}

Let $\chi_k$ be the character of the representation $\wedge^k R$. By orthogonality of characters, we have 
\begin{equation}
\label{numero1} T(A^*, k, a^*) = \frac{1}{|A|} \sum_{a} a^{*-1}(a) \chi_k(a).
\end{equation}

In other words, the function $T(A^*, k, a^*)$ on $A^*$ is the Fourier transform of the function $\chi_k$ on $A$.

\begin{proposition} \label{propositionA2} The value of the character $\chi_k$ on $a$ in $A$ depends
 only on the order $d$ of $a$. It is given by:
\[
\chi_k(a) = \begin{cases}
(-1)^{(d-1)(k/d)} \binom{|A|/d}{k/d} & \text{if } d \mid k; \\
0 & \text{otherwise}.
\end{cases}
\]
\end{proposition}

\begin{proof}
For $B$ a $k$-element subset of $A$, the exterior product of the $\delta_\beta$ ($\beta \in B$) 
depends, but only up to sign, on choosing an ordering of $B$. Let $L_B$ be the line spanned by 
this exterior product. It determines $B$, and the exterior power $\wedge^k R$ is the direct sum 
of the lines $L_B$. If $L_B$ is stable by $a$, so is $B$, and $B$ is the union of orbits of the cyclic 
group $\langle a \rangle$ generated by $a$. As $\langle a \rangle$ acts freely, it follows that $d$ 
divides $k$ and that $a$ acts on $L_B$ by $(-1)^{(k/d)(d-1)}$. 
Such $B$ are the
inverse images of the $k/d$-element subsets of $A / \langle a \rangle$.  
The formula for the trace $\chi_{k}(a)$ of the action of $a$ on $\wedge^{k}R$ follows. 
\end{proof} 

The abelian group $A$ is the direct sum of cyclic groups $C_i$ of order $m_i$; 
$(1 \leq i \leq c)$, with $m_{i+1} \mid m_i$. 
Define $n := m_1$ to be the maximum order of elements of $A$. 

For $d$ a divisor of $n$, let $A_d$ be the group $\ker(d: A \to A)$ of elements of $A$ of order dividing $d$. 
The set $\text{Div}(n)$ of divisors of $n$, ordered by divisibility, is a lattice, 
with $\inf = \gcd$ and $\sup = \text{lcm}$. It is isomorphic, by $d \mapsto A_d$, 
to the set of the subgroups $A_d$ of $A$, ordered by inclusion. In the latter lattice, 
$\inf$ is intersection and $\sup$ is ``subgroup generated by".  

Pontryagin duality is an exact functor. It follows that the orthogonal in $A^*$ of $A_d$
 is $d A^* := \text{Im}(d: A^* \to A^*)$. Indeed, $\ker(d: A \to A)^*$ is $\text{coker}(d: A^* \to A^*)$ 
 and $\text{Im}(d: A^* \to A^*)$ is the kernel of the map from $A^*$ to $\ker(d: A \to A)^*$. 
 The subgroups $d A^*$, ordered by inclusion, form again a lattice, with $\inf =$ intersection, 
 $\sup =$ group generated by. For $a^* \in A^*$, let $d^*(a^*)$ be the largest $d\mid n$ 
 such that $a^*$ is in $d A^*$.  

\begin{proposition} \label{propositionA3} The Fourier transform maps the space of functions $\chi$ on $A$ such that $\chi(a)$ 
depends only on the order of $a$ onto the space of functions $T$ on $A^*$ such that $T(a^*)$ 
depends only on $d^*(a^*)$.  
\end{proposition}

\begin{proof} Because the set of subgroups $A_d$ of $A$ (resp.\ $d A^*$ of $A^*$) 
is stable by intersections, exclusion-inclusion shows that the spaces of functions 
considered on $A$ (resp.\ on $A^*$) admit as basis the characteristic functions of the 
$A_d$ (resp.\ $d A^*$) for $d \mid n$. It remains to observe that the Fourier transform 
maps the characteristic function of $A_d$ to $\frac{|A_d|}{|A|}$ times the characteristic 
function of the orthogonal $d A^*$ of $A_d$.
\end{proof}

Let us compute the matrix of Fourier transform, restricted to spaces considered in 
Proposition \ref{propositionA3}, when using as basis the $e_d$ (resp.\ $e_d ^*$), where $e_d$ 
(resp.\ $e_d ^*$) is the characteristic function of the set of elements $a$ of $A$ 
(resp.\ $a^*$ of $A^*$) of order $d$ (resp.\ for which $d^*(a^*) = d$).

For each prime $p$ dividing $n$, let $A[p]$ be the $p$-Sylow subgroup of $A$, 
and $n[p]$ be the corresponding $n$. We have $A = \bigoplus A[p]$, $n = \prod n[p]$, 
the lattice of divisors of $n$ is the product of the totally ordered sets of the divisors of
the $n[p]$, the space of the functions on $A$ (resp.\ $A^*$) depending only on the order
 (resp.\ $d^*(a^*)$) is the tensor product of the corresponding spaces on the $A[p]$ 
 (resp.\  $A^*[p] = A[p]^*$), the bases defined above are the tensor product of the 
 corresponding bases for the $A[p]$ (resp.\ $A[p]^*$) 
and the Fourier transform for $A$ is the tensor product of the Fourier transforms 
for the $A[p]$. It follows that the matrix we want to compute is the Kronecker product 
of the similar matrices for the $A[p]$. 

This reduces us to the case where $A$ is an abelian $p$-group, and $n$ a power $p^m$ of $p$. 
Let us write $e_a$ (resp.\ $e_a ^*$) for what was called $e_{p^a}$ 
(resp.\ $e^*_{p^a}$) $(0 \leq a \leq m)$. We have
\begin{equation}
\begin{aligned}
e_a &= \text{characteristic function of } A_{p^a} \\
&-{} \text{characteristic function of } A_{p^{a-1}}.
\end{aligned}
\end{equation}

Hence, the Fourier transform of $e_a$ is:
\begin{equation}
\begin{aligned}
&\frac{|A_{p^a}|}{|A|} \text{ characteristic function of } p^a A^* \\
&\quad - \frac{|A_{p^{a-1}}|}{|A|} \text{ characteristic function of } p^{a-1} A^* \\
&= \frac{1}{|A|} \left[ - |A_{p^{a-1}}| e_{a-1}^* + \left( |A_{p^a}| - |A_{p^{a-1}}| \right) \sum_{b=a}^{m} e^*_b \right].
\end{aligned}
\end{equation}

When $A$ is cyclic of prime power order $p^m$, one has $p^a A^* = A_{p^{m-a}} ^*$ 
and one recovers the formula of Conjecture 1.

\subsection*{Variant}

A multiset $B$ of a set $I$  is a function (the \textit{multiplicity}) $\varphi: I \to \mathbb{N}$. 
The \textit{cardinality} of $B$ is the sum of the multiplicities. 
If $e_i$ ($i \in I$) is a basis of $R$, the products $\prod_{i \in B} e_i$, for $B$ of cardinality $k$, 
are a basis of $\text{Sym}^k(R)$.

Define $T^+(A, k, a)$ to be the number of multisets of $A$ of cardinality $k$ summing up to $a$.

Repeating the arguments proving Proposition \ref{propositionA1}, one finds that the function $T^+(A^*, k, a^*)$ 
on $A^*$ is the Fourier transform of the character $\chi_k^+$ of the representation $\text{Sym}^k(R)$, 
and that for $a$ in $A$ of order $d$,
$$\chi_k^+(a) = \begin{cases}
\binom{|A|/d + k/d-1}{k/d} & \text{if } d | k;\\
0 & \text{otherwise.}
\end{cases}$$
The binomial coefficient is the number of multisets of cardinality $k/d$ of $A/\langle a\rangle$. 
It can be rewritten as 
$$(-1)^{k/d}\binom{-|A|/d}{k/d}.$$

\subsection*{Generalization}

We have used the functors of exterior and symmetric powers. 
Other polynomial functors: (vector spaces) $\longrightarrow$ (vector spaces) could be used.

Let $S_N$ be the algebra of symmetric polynomials with integer coefficients in $N$ 
variables $X_1, \ldots, X_N$. A polynomial representation of the algebraic group $GL(N)$
defines an element of $S_N$: the restriction of its character to the diagonal matrices 
$\text{diag}(X_1, \ldots, X_N)$. This construction gives an isomorphism from the
Grothendieck group of the category of polynomial representations of $GL(N)$ to $S_N$. 
The elements of this Grothendieck group will be called virtual representations.

To the exterior power $\wedge^k$ (resp.\ symmetric power $\text{Sym}^k$) of the standard representation 
correspond the elementary symmetric polynomial $e_k$ (resp.\ the complete symmetric polynomial $h_k$).

Take $N = |A|$, and let $\rho$ be a polynomial representation of $GL(N)$, with corresponding 
symmetric polynomial $P$. The representation $\rho$ defines a functor from $N$-dimensional
 vector spaces to vector spaces, noted $V \mapsto \rho(V)$. We will use the representation 
 $\rho(R)$ of $A$. In the basis $a^*$ of $R$, $A$ acts by diagonal matrices, so that we can 
 read the character of the representation $\rho(R)$ of $A$ from $P$, as follows. Let us assign 
 to each variable $X_i$ a different element of $A^*$. This defines an homomorphism from 
 $S_N$ to the group algebra of $A^*$, and this homomorphism does not
depend on the assignment used. Let $\sum n_P(a^*) a^*$ be the image of $P$. The function 
$$\chi_P: a \mapsto \sum n_P(a^*) a^*(a)$$
is the character of $\rho(R)$, and $n_P(a^*)$ is the function on $A^*$ Fourier transform of this character.

The construction, being additive in $\rho$, extends to virtual representations of $GL(N)$.

\begin{example}
(i) Take for $P$ the monomial symmetric function
$$m_\lambda = X_1^{\lambda_1} X_2^{\lambda_2} \ldots X_N^{\lambda_N} + \text{ distinct conjugates}.$$
Let $B_\lambda$ be the set of multisets of $A^*$ with multiset of multiplicities $\lambda$. 
Then $n_P(a^*)$ is the number of elements of $B_\lambda$ summing up to $a^*$.\\
(ii) For $P = e_k$ (resp.\ $h_k$), $n_P(a^*)$ is $T(A^*, k, a^*)$ (resp.\ $T^+(A^*, k, a^*)$).
\end{example}

If $a$ in $A$ is of order $d$, the eigenvalues of $a$ acting on $R$ are the $d^{th}$ roots of 1, 
each having multiplicity $|A|/d$. If one gives the $X_i$ those values (with $|A|/d$ variables 
assigned to each $d^{th}$ root of 1), $\chi_P(a)$ is the corresponding value of $P$.

\begin{corollary} 
(i) $\chi_P(a)$ depends only on the order of $a$.\\
(ii) $n_P(a^*)$ depends only on $d^*(a^*)$.\\
(iii) If $P$ is homogeneous of degree $k$,
$\chi_P(a)$ vanishes if the order $d$ of $a$ does not divide $k$ 
(and hence the $\gcd(n, k)$). 
Dually, $n_P(a^*)$ depends only on the image of $a^*$ in $A^*/\gcd(n,k)A^*$. 
\end{corollary}

\begin{proof}
(i) is clear and (ii) follows from Proposition \ref{propositionA3}.
The system of eigenvalues for the action of $a$ on $R$ is stable 
by multiplication by any $d^{th}$ root of unity $\zeta$. It follows that for $P$ 
homogeneous of degree $k$, the corresponding value $\chi_P(a)$ 
satisfies $\chi_P(a) = \zeta^k \chi_P(a)$, proving the first assertion of (iii). 
The second is its Fourier transform.
\end{proof}

Applying the construction above to the elementary symmetric function 
$e_k$, one gets another proof of Proposition \ref{propositionA2}. Indeed, for $P = e_k$, 
$n_P(a)$, for $a$ of order $d$, is $(-1)^k$ times the coefficient of $T^{|A|-k}$ 
in the product over the eigenvalues $\alpha$ of the $T - \alpha$. This product is 
$$(T^d-1)^{|A|/d}$$
which one expands by the binomial formula.

\noindent 2020 {\it Mathematics Subject Classification}: Primary 11B30; Secondary 05A18, 81T18, 20C15. 

\noindent \emph{Keywords:} Feynman diagram, Lyndon word, necklace, partition, Ramanujan sum, virtual representation. 

\noindent(Concerned with sequences 
\seqnum{A001840}, \seqnum{A006918}, \seqnum{A007997}, \seqnum{A008610}, 
\seqnum{A008646}, \seqnum{A011795}, \seqnum{A011796}, \seqnum{A011797}, 
\seqnum{A031164}, \seqnum{A032191}, \seqnum{A032192}, \seqnum{A032193}, 
\seqnum{A032194}, \seqnum{A032801}, \seqnum{A051168}, \seqnum{A054535}, 
\seqnum{A263318}, \seqnum{A267632}, \seqnum{A379920}, \seqnum{A381289}, 
\seqnum{A381290}, \seqnum{A381291}, \seqnum{A381350}, and \seqnum{A381351}.)  
 
\end{document}